\documentclass[letterpaper, 10pt, conference]{IEEEconf}
\IEEEoverridecommandlockouts 

\usepackage{amsmath,amssymb}
\usepackage{epsfig,subfigure,url,warmread}
\usepackage[all,import]{xy}

\newcommand{\norm}[1]{\ensuremath{\left\| #1 \right\|}}
\newcommand{\abs}[1]{\ensuremath{\left| #1 \right|}}
\newcommand{\bracket}[1]{\ensuremath{\left[ #1 \right]}}
\newcommand{\braces}[1]{\ensuremath{\left\{ #1 \right\}}}
\newcommand{\parenth}[1]{\ensuremath{\left( #1 \right)}}
\newcommand{\refeqn}[1]{(\ref{eqn:#1})}
\newcommand{\reffig}[1]{Fig. \ref{fig:#1}}
\newcommand{\tr}[1]{\mbox{tr}\ensuremath{\negthickspace\bracket{#1}}}
\newcommand{\deriv}[2]{\ensuremath{\frac{\partial #1}{\partial #2}}}
\newcommand{\SO}{\ensuremath{\mathrm{SO(3)}}}
\newcommand{\T}{\ensuremath{\mathrm{T}}}
\newcommand{\so}{\ensuremath{\mathfrak{so}(3)}}

\renewcommand{\Re}{\ensuremath{\mathbb{R}}}
\renewcommand{\S}{\ensuremath{\mathbb{S}}}

\title{\LARGE \bf
Global Attitude Estimation using Single Direction Measurements}

\author{Taeyoung Lee\authorrefmark{1}\authorrefmark{2}, Melvin Leok\authorrefmark{1}, N. Harris McClamroch\authorrefmark{2}, and Amit Sanyal%
\thanks{Taeyoung Lee, Aerospace Engineering, University of Michigan, Ann Arbor, MI 48109 {\tt tylee@umich.edu}}%
\thanks{Melvin Leok, Mathematics, Purdue University, West Lafayette, IN 47907 {\tt mleok@math.purdue.edu}}%
\thanks{N. Harris McClamroch, Aerospace Engineering, University of Michigan, Ann Arbor, MI 48109 {\tt nhm@umich.edu}}%
\thanks{Amit Sanyal, Mechanical and Aerospace Engineering, Arizona State University, Tempe, AZ
85287 {\tt sanyal@asu.edu}}%
\thanks{\textsuperscript{\footnotesize\ensuremath{*}}This research has been supported in part by NSF under grant DMS-0504747, and by a grant from the Rackham Graduate School, University of Michigan.}
\thanks{\textsuperscript{\footnotesize\ensuremath{\dagger}}This research has been supported in part by NSF under grant ECS-0244977.}
}

\begin{document}
\maketitle \thispagestyle{empty} \pagestyle{empty}

\begin{abstract}
A deterministic attitude estimator for a rigid body under an
attitude dependent potential is studied. This estimator requires
only a single direction measurement to a known reference point at
each measurement instant. The measurement cannot completely
determine the attitude, but an attitude estimation scheme based on
this measurement is developed; a feasible set compatible with the
measurement is described and it is combined with an attitude
dynamics model to obtain an attitude estimate. The attitude is
globally represented by a rotation matrix, and the uncertainties are
described by ellipsoidal sets. A numerical example for a spacecraft
in a circular orbit is presented.
\end{abstract}

\section{Introduction}

The attitude of a rigid body is defined by the orientation of a
body-fixed frame with respect to a reference frame, and the attitude
is represented by a rotation matrix that is a $3\times 3$ orthogonal
matrix with determinant of 1, which transforms a representation of a
vector in a body-fixed frame into one represented in the reference
frame. Rotation matrices have a group structure denoted by $\SO$. In
spacecraft applications, the attitude is usually determined by using
a set of direction measurements. The directions to objects such as
the sun, stars, and geomagnetic fields, assumed to be known in the
reference frame, are measured in the body-fixed frame in order to
determine the rotation matrix.

Attitude determination using multiple direction measurements with
least squares estimation is known as Wahba's
problem~\cite{jo:wahba}. The original solution of Wahba's problem is
given in~\cite{jo:solwahba}, and solutions are expressed in terms of
quaternions~\cite{jo:ShOh1981}, and in terms of a rotation
matrix~\cite{pro:san2006}. Based on these attitude determination
schemes, attitude estimation problems are studied in
\cite{Shu.JAS90,Pi.JGCD00} and \cite{CDC06est}. The attitude
determination/estimation procedures using Wahba's problem
formulation require at least two different direction measurements at
each measurement instant. This places a stringent constraint on the
dynamic estimation of spacecraft attitude.

A single direction measurement provides some information about the
attitude; it is  guaranteed that the rotation matrix lies in a one
dimensional subgroup of the three dimensional special orthogonal
group \SO, which is diffeomorphic to the one-sphere $\S^1$. The
attitude is not completely determined at a single measurement
instant. If the process is coupled with an attitude dynamics model,
an attitude estimation scheme can be developed using single
direction measurements. An attitude determination scheme using
single direction measurements is studied in~\cite{Shu.JAS04}, but
this approach requires an additional arc length measurement.

Most existing attitude estimation schemes use generalized coordinate
representations of the attitude. As is well known, minimal
coordinate representations of the rotation group, such as Euler
angles, lead to singularities. Non-minimal coordinate
representations, like the quaternions, have their own associated
problems. Besides the extra constraint of unit norm that one needs
to impose on the quaternion, the quaternion representation, which is
diffeomorphic to $\mathrm{SU(2)}$, double covers $\SO$. So, it has
an inevitable ambiguity in expressing the attitude.

A stochastic state estimator requires probabilistic models for the
state uncertainty and the noise. However, statistical properties of
the uncertainty and the noise are often not available. An
alternative deterministic approach is to specify bounds on the
uncertainty and the measurement noise without an assumption on their
distribution. Noise bounds are available in many cases, and
deterministic estimation is robust to the noise
distribution~\cite{jo:TheSkaSou.IEEETSP94}. An efficient but
flexible way to describe the bounds is using ellipsoidal sets,
referred to as uncertainty ellipsoids. The deterministic estimation
process is based on set theory results developed
in~\cite{jo:Sc1968}; optimal deterministic estimation problems using
the uncertainty ellipsoids are studied in~\cite{jo:MaNo1996} and
\cite{jo:DuWaPo2001}.

In this paper, a deterministic attitude estimator which requires a
single direction measurement at each measurement instant is
presented. A feasible set in $\SO$ that is compatible with the
measurement is represented by Lie algebra elements and the
exponential map. It is compared with the attitude dynamics model to
obtain an updated estimate. The estimation scheme presented in this
paper has the following distinctive features: the estimator requires
only a single direction measurement at each measurement instant, the
attitude is represented by a rotation matrix without any local
parameterization, and the deterministic estimator is distinguished
from a Kalman or extended Kalman filter. 

This paper is organized as follows. The attitude dynamics and
uncertainty model are given in Section II. The attitude
determination scheme and the attitude estimation scheme using single
direction measurements are presented in Section III and IV, which is
followed by a numerical example in Section V.

\section{Attitude dynamics and Uncertainty model}

\subsection{Equations of motion}
We consider estimation of the attitude dynamics of a rigid body in
the presence of an attitude dependent potential,
$U(\cdot):\SO\mapsto\Re$, $R\in\SO$. Systems that can be so modeled
include a free rigid body, spacecraft on a circular orbit with
gravity gradient effects~\cite{pro:acc06}, or a 3D
pendulum~\cite{pro:cca05}. The continuous equations of motion are
\begin{gather}
J\dot\Omega + \Omega\times J\Omega = M,\\
\dot{R} = R S(\Omega),\label{eqn:Rdot}
\end{gather}
where $J\in\Re^{3\times 3}$ is the moment of inertia matrix of the
rigid body, $\Omega\in\Re^3$ is the angular velocity of the body
expressed in the body-fixed frame, and $S(\cdot):\Re^3\mapsto \so$
is a skew mapping defined by $S(x)y=x\times y$ for all
$x,y\in\Re^3$. The vector $M\in\Re^3$ is the moment due to the
potential, determined by $S(M)=\deriv{U}{R}^TR-R^T\deriv{U}{R}$, or
more explicitly,
\begin{gather}
M=r_1\times v_{r_1} + r_2\times v_{r_2} +r_3\times v_{r_3},
\end{gather}
where $r_i,v_{r_i}\in\Re^{1\times 3}$ are the $i$th row vectors of
$R$ and $\deriv{U}{R}$, respectively.

General numerical integration methods like the popular Runge-Kutta
schemes, typically preserve neither first integrals nor the
characteristics of the configuration space, $\SO$. In particular,
the orthogonal structure of the rotation matrices is not preserved
numerically. It is often proposed to parameterize \refeqn{Rdot} by
Euler angles or quaternions instead of integrating \refeqn{Rdot}
directly. However, Euler angles yield only local representations of
the attitude and they have singularities. Unit quaternions do not
exhibit singularities, but they have the manifold structure of the
three-sphere $\S^3$, and double cover $\SO$. Consequently, the unit
quaternion representing the attitude is inevitably ambiguous. In
addition, general numerical integration methods do not preserve the
unit length constraint. Therefore, quaternions have the same
numerical drift problem as rotation matrices.

Lie group variational integrators preserve the group structure
without the use of local charts, reprojection, or constraints, they
are symplectic and momentum preserving, and they exhibit good energy
behavior for an exponentially long time period. The following Lie
group variational integrator for the attitude dynamics of a rigid
body is presented in~\cite{pro:cca05}:
\begin{gather}
h S(J\Omega_k+\frac{h}{2} M_k) = F_k J_d - J_dF_k^T,\label{eqn:findf0}\\
R_{k+1} = R_k F_k,\label{eqn:updateR0}\\
J\Omega_{k+1} = F_k^T J\Omega_k +\frac{h}{2} F_k^T M_k
+\frac{h}{2}M_{k+1},\label{eqn:updatew0}
\end{gather}
where $J_d\in\Re^{3\times 3}$ is a nonstandard moment of inertia
matrix defined by $J_d=\frac{1}{2}\mathrm{tr}\!\bracket{J}I_{3\times
3}-J$, and $F_k\in\SO$ is the relative attitude between integration
steps. The constant $h\in\Re$ is the integration step size, and the
subscript $k$ denotes the $k$th integration step. This integrator
yields a map $(R_k,\Omega_k)\mapsto(R_{k+1},\Omega_{k+1})$ by
solving \refeqn{findf0} to obtain $F_k\in\SO$ and substituting it
into \refeqn{updateR0} and \refeqn{updatew0} to obtain $R_{k+1}$ and
$\Omega_{k+1}$. The only implicit part is \refeqn{findf0}. The
actual computation of $F_k$ is done in the Lie algebra $\so$ of
dimension 3, and the rotation matrices are updated by
multiplication. So this approach is distinguished from integration
of the kinematics equation \refeqn{Rdot}, and there is no excessive
computational burden. We use these discrete equations of motion to
propagate the attitude dynamics between measurements during the
estimation process.

\subsection{Uncertainty Ellipsoid}
We describe uncertainties of the attitude dynamics by using
ellipsoidal sets referred to as uncertainty ellipsoids. An
uncertainty ellipsoid in $\Re^n$ is defined as
\begin{align}
    \mathcal{E}_{\Re^n}(\hat x,P)=\braces{x\in\Re^n \,\Big|\,
    (x-\hat{x})^T P^{-1}(x-\hat{x})\leq
    1},
\end{align}
where $\hat{x}\in\Re^n$, and $P\in\Re^{n\times n}$ is a symmetric
positive definite matrix. We call $\hat x$ the center of the
uncertainty ellipsoid, and $P$ is the uncertainty matrix that
determines the size and the shape of the uncertainty ellipsoid. The
size of an uncertainty ellipsoid is measured by
$\mathrm{tr}\!\bracket{P}$ which is the sum of the squares of the
semi principal axes of the ellipsoid.

The attitude dynamics evolves on the 6 dimensional tangent bundle,
$\T\SO$. We identify $\T\SO$ with $\SO\times\so$ by left
trivialization, and we identify $\so$ with $\Re^3$ by the
isomorphism $S(\cdot)$. The uncertainty ellipsoid centered at
$(\hat{R},\hat\Omega)\in\T\SO$ is induced from an uncertainty
ellipsoid in $\Re^6$;
\begin{align}
    \mathcal{E}(\hat{R},\hat{\Omega},P)
        & = \braces{R\in\SO,\,\Omega\in\Re^3 \,\Big|\,
    \begin{bmatrix}\zeta\\\delta\Omega\end{bmatrix}\in\mathcal{E}_{\Re^6}(0_{6},P)},\label{eqn:ueso}
\end{align}
where $S(\zeta)=\mathrm{logm} \parenth{\hat{R}^T R}\in\so$,
$\delta\Omega=\Omega-\hat{\Omega}\in\Re^3$, and $P\in\Re^{6\times
6}$ is a symmetric positive definite matrix. An element $(R,\Omega)
\in\mathcal{E}(\hat{R},\hat{\Omega},P)$ can be written as
\begin{align*}
    R = \hat{R} \exp{ S(\zeta)},\quad
    \Omega = \hat{\Omega} + \delta \Omega,
\end{align*}
for some $x=\bracket{\zeta ;\,\delta\Omega}\in\Re^6$ satisfying
$x^TP^{-1}x\leq 1$.

We assume that the initial conditions are bounded by a prescribed
uncertainty ellipsoid
\begin{gather}
    (R_0,\Omega_0)\in\mathcal{E}(\hat{R}_0,\hat{\Omega}_0,P_0),\label{eqn:P0}
\end{gather}
where $P_0\in\Re^{6\times 6}$ is a symmetric positive definite
matrix that defines the shape and the size of the uncertainty
ellipsoid.

\section{Attitude Determination with a Single Direction
Measurement}\label{sec:mea}

In the attitude determination problem, we measure directions to
points in the reference frame. We assume that the directions to
these points are known in the reference frame. This either requires
that the points are located far away from the spacecraft or the
relative location of the spacecraft is known exactly. The
directional sensor is fixed in the body-fixed frame, and the
measurements are representations of the direction vectors in the
body-fixed frame. The representations in the body-fixed frame are
transformed into those in the reference frame by multiplication with
the rotation matrix that defines the attitude of the rigid body.

\subsection{Exact measurement}
Let the direction to a known point in the reference frame be
$e\in\S^2$, and let the corresponding vector represented in the
body-fixed frame be $b\in\S^2$. We first assume that the direction
measurement has no error, so the direction $b$ is exact. Since we
only measure a direction to a point, we normalize $e$ and $b$ so
that they have unit lengths. The
 vectors $e$ and $b$ are different representations of the same vector
 from the spacecraft to the known point, and they are
 related by a rotation matrix $R\in\SO$
that defines the attitude of the rigid body
\begin{align}
e = R b.\label{eqn:eRb}
\end{align}

This equation provides a two-dimensional constraint on the
three-dimensional rotation matrix. Consequently, a single direction
measurement does not completely determine the attitude. This
corresponds to the fact that if we rotate the rigid body about the
direction $e$ in the reference frame, then the measured direction
$b$ is not changed. The rotation matrix has one-dimensional
uncertainty represented by any rotation about the direction $e$ in
the reference frame, or equivalently, any rotation about the
direction $b$ in the body-fixed frame.

Suppose that $R^\circ\in\SO$ is a particular rotation matrix
satisfying \refeqn{eRb}. This rotation matrix can be represented in
several ways. For example, if $b$ and $e$ are not co-linear,
\begin{align}\label{eqn:Rcirc}
    R^\circ=\exp\!\bracket{\cos^{-1}(b^Te)\,\, S\!\parenth{\frac{b\times
    e}{\norm{b\times e}}}}\exp S(\theta^\circ b),
\end{align}
where the constant $\theta^\circ \in \S^1$ can be arbitrarily
chosen. The rotation matrix that represents the attitude of the
rigid body can be written in terms of $R^\circ$ as
\begin{align}\label{eqn:Ract}
    R= R^\circ \exp\bracket{\theta S(b)}
\end{align}
for a $\theta\in\S^1$.

In summary, if the single direction to a known point is measured
exactly, the rotation matrix lies in the following one dimensional
subgroup of $\SO$:
\begin{align}
    R\in\braces{R^\circ \exp\bracket{\theta S(b)}\bigg|
    \theta\in\S^1}.
\end{align}

\subsection{Measurement error}
We now consider the effects of small measurement errors. Let $\tilde
b\in\S^2$ be the measured direction of the direction $b$. Since we
only measure directions, we normalize $b$ and $\tilde b$ so that
they have unit lengths. Therefore it is inappropriate to express the
measurement error by a vector difference. The measurement error is
modeled by rotation of the measured direction
\begin{align}
{b}& = \exp\bracket{S(\nu)} \tilde b,\nonumber\\
& \simeq \tilde b + S(\nu)\tilde b,\label{eqn:bi}
\end{align}
where $\nu\in\Re^3$ is the Euler axis of rotation from $\tilde b$ to
$b$, and $\norm{\nu}$ is the corresponding rotation angle error in
radians. The measurement error is bounded by an uncertainty
ellipsoid
\begin{align}
    \nu \in \mathcal{E}_{\Re^3}(0_3,S)
\end{align}
for a symmetric positive definite matrix $S\in\Re^{3\times 3}$. The
magnitude of the measurement error is assumed to be small.

Let $\tilde R^\circ\in\SO$ be a rotation matrix obtained by
\refeqn{Rcirc} for the measured direction $\tilde b$. 
We express the difference between $R^\circ$ and $\tilde R^\circ$
using the exponential map:
\begin{align}\label{eqn:Rcircerr}
    R^\circ=\tilde R^\circ \exp\bracket{S(\zeta^\circ)}
\end{align}
for some $\zeta^\circ\in\Re^3$. Since we make the small measurement
error assumption, the norm of the vector $\zeta^\circ$ is considered
to be much smaller than $\pi$, i.e. $\norm{\zeta^\circ}\ll \pi$.
Since $e=R^\circ b =\tilde R^\circ\tilde b$, we obtain
\begin{align*}
    \tilde R^\circ \tilde b & = R^\circ b,\\
    & = \tilde R^\circ \exp\bracket{S(\zeta^\circ)}\braces{I_{3\times 3}+S(\nu)}\tilde
    b,\\
    & \simeq \tilde R^\circ \braces{I_{3\times 3}+S(\zeta^\circ+\nu)}\tilde
    b,
\end{align*}
Thus we have $S(\zeta^\circ+\nu)\tilde b=0$, which is equivalent to
\begin{align}\label{eqn:zetacirc}
    \zeta^\circ = c\tilde b - \nu
\end{align}
for any constant $c\in\Re$. Since $\norm{\zeta^\circ}\leq \abs{c}
+\norm{\nu}\ll \pi$, the constant $c$ is smaller than $\pi$, i.e.
$\abs{c}\ll\pi$.

Substituting \refeqn{bi}, \refeqn{Rcircerr}, and \refeqn{zetacirc}
into \refeqn{Ract}, we obtain
\begin{align}\label{eqn:Ractm}
    R = \tilde R^\circ
    \exp\!\bracket{S(c\tilde b - \nu)}\exp\!\bracket{\theta S\parenth{(I_{3\times
    3}+S(\nu))\tilde b}}
\end{align}
for constants $c$ and $\theta\in\S^1$.

In summary, if the single direction measurement has a small error
represented by \refeqn{bi}, then the attitude of the rigid body can
be written in terms of the measured direction $\tilde b$ and the
measurement error $\nu$ as \refeqn{Ractm}. This expression includes
the uncertainty caused by the measurement error as well as the
uncertainty due to the single direction measurement assumption. The
constant $\theta^\circ\in\S^1$ to determine $R^\circ$ and $\tilde
R^\circ$ is specified by the following estimation procedure.

\section{Attitude Estimation with a Single Direction
Measurement}\label{sec:est} The deterministic estimation scheme
using uncertainty ellipsoids is introduced first. A deterministic
estimator for the attitude and the angular velocity of a rigid body
is developed by using the preceding attitude determination scheme.

The subscript $k$ denotes the $k$th discrete index. The superscript
$f$ denotes the variables related to the flow update, and the
superscript $m$ denotes the variables related to the measurement
update. $\tilde\cdot$ denotes a measured variable, and $\hat\cdot$
denotes an estimated variable.

\subsection{Deterministic estimation}

We use deterministic bounded estimation using ellipsoidal sets,
referred to as uncertainty ellipsoids, to describe the uncertainty
and measurement noise. The estimation process has three steps
similar to those in the Kalman filter: prediction, measurement, and
filtering steps. We assume that the initial condition lies in a
prescribed uncertainty ellipsoid, which is propagated in time using
the equations of motion. This defines a prediction step. The
measurement error bound is described by a measurement uncertainty
ellipsoid. Then we can guarantee that the state lies in the
intersection of the predicted uncertainty ellipsoid and the measured
uncertainty ellipsoid. The intersection of the two ellipsoids is an
irregular shape, which is not efficient to compute and store.
Instead we find a minimal ellipsoid that contains this intersection.
This procedure is repeated whenever new measurements are available.

This deterministic estimation procedure is illustrated in
\reffig{ue}. The left figure shows time evolution of an uncertainty
ellipsoid, and the right figure shows a cross section at a fixed
measurement instant. At the $k$th time step, the state is bounded by
an uncertainty ellipsoid centered at $\hat{x}_k$. This initial
ellipsoid is propagated through time. Suppose that the state is
measured next at the $(k+l)$th time step. In this single direction
measurement estimation, the measurement ellipsoid degenerates to a
strip. At this instant, the actual state lies in the intersection.
In the estimation process, we find a new ellipsoid that contains
this intersection, as shown in the right figure. The center of the
new ellipsoid, $\hat{x}_{k+l}$ gives a point estimate of the state
at time step $(k+l)$, and the magnitude of the new uncertainty
ellipsoid measures the estimation accuracy. The deterministic
estimates are optimal in the sense that the sizes of the ellipsoids
are minimized.

\renewcommand{\xyWARMinclude}[1]{\includegraphics[width=0.30\textwidth]{#1}}
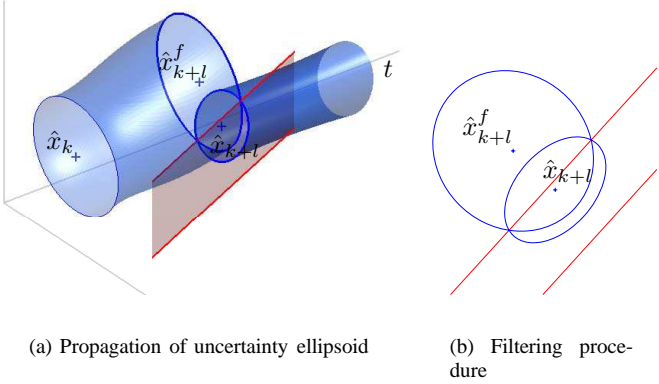
\begin{figure}[t]
    \centerline{\subfigure[Propagation of uncertainty ellipsoid]{
    $$\begin{xy}
    \xyWARMprocessEPS{tubedeg}{eps}
    \xyMarkedImport{}
    \xyMarkedMathPoints{1-5}
    \end{xy}$$}
    \renewcommand{\xyWARMinclude}[1]{\includegraphics[width=0.17\textwidth]{#1}}
    \subfigure[Filtering procedure]{
        $$\begin{xy}
    \xyWARMprocessEPS{tubesecdeg}{eps}
    \xyMarkedImport{}
    \xyMarkedMathPoints{1-3}
    \end{xy}$$}}
    \caption{Uncertainty ellipsoids}\label{fig:ue}
\end{figure}

\subsection{Flow update}
Suppose that the attitude and the angular momentum at the $k$th step
lie in a given uncertainty ellipsoid:
\begin{align*}
    (R_k,\Omega_k)\in\mathcal{E}(\hat{R}_k,\hat{\Omega}_k,P_k),
\end{align*}
and a new measurement is taken at the $(k+l)$th time step.

 The flow
update finds the center and the uncertainty matrix that define the
uncertainty ellipsoid at the $(k+l)$th step using the given
uncertainty ellipsoid at the $k$th step. Since the attitude dynamics
of a rigid body is nonlinear, the admissible boundary of the state
at the $(k+l)$th step is not an ellipsoid in general. We assume that
the given uncertainty ellipsoid at the $k$th step is sufficiently
small that attitudes and angular velocities in the uncertainty
ellipsoids can be approximated using the linearized equations of
motion. Then we can guarantee that the uncertainty set at the
$(k+l)$th step is an ellipsoid, and we can compute its center and
its uncertainty matrix at the $(k+l)$th step separately.

 \textit{Center:} For the given center at step $k$, $(\hat{R}_k,\hat{\Omega}_k)$, the center of the uncertainty
ellipsoid at step $(k+l)$ is
$(\hat{R}_{k+l}^{f},\hat{\Omega}_{k+l}^{f})$ obtained using the
discrete equations of motion, \refeqn{findf0}, \refeqn{updateR0},
and \refeqn{updatew0}:
\begin{gather}
h S(J\hat{\Omega}_k+\frac{h}{2} \hat{M}_k) = \hat{F}_k J_d - J_d
\hat{F}_k^T,\label{eqn:findf}\\
\hat{R}_{k+1}^{f} = \hat{R}_k \hat{F}_k,\label{eqn:updateR}\\
J\hat{\Omega}_{k+1}^{f} = \hat{F_k}^T
\hat{\Omega}_k+\frac{h}{2}\hat{F_k}^T \hat{M_k} +\frac{h}{2}
\hat{M}_{k+1}.\label{eqn:updatePi}
\end{gather}
This integrator yields a map $(\hat R_k,\hat\Omega_k)\mapsto(\hat
R^f_{k+1},\hat\Omega_{k+1}^f)$, and this process is repeatedly
applied to find the center at the $(k+l)$th step, $(\hat
R^f_{k+l},\hat\Omega_{k+l}^f)$.

\textit{Uncertainty matrix:} We assume that an uncertainty ellipsoid
contains small perturbations from the center of the uncertainty
ellipsoid. Then the uncertainty matrix is propagated by using the
linearized flow of the discrete equations of motion. At the
$(k+1)$th step, the uncertainty ellipsoid is represented by
perturbations from the center $(\hat R^f_{k+1},\hat\Omega_{k+1}^f)$
as
\begin{align*}
    R_{k+1}&=\hat{R}_{k+1}^{f} \exp{S(\zeta_{k+1}^{f})},\\
    \Omega_{k+1}&=\hat{\Omega}_{k+1}^{f}+\delta\Omega_{k+1}^{f},
\end{align*}
for some $\zeta_{k+1}^{f},\delta\Omega_{k+1}^{f}\in\Re^3$. The
uncertainty matrix at the $(k+1)$th step is obtained by finding a
bound on $\zeta_{k+1}^{f},\delta\Omega_{k+1}^{f}\in\Re^3$. Assume
that the uncertainty ellipsoid at the $k$th step is sufficiently
small. Then, $\zeta_{k+1}^{f},\delta\Omega_{k+1}^{f}$ are
represented by the following linear equations using the results
presented in~\cite{pro:acc06}
\begin{align*}
x_{k+1}^{f} & = A_k^f x_k,
\end{align*}
where $x_k=[\zeta_k;\delta\Omega_k]\in\Re^6$, and
$A_k^f\in\Re^{6\times 6}$ can be suitably defined. Since
$(R_k,\Omega_k)\in\mathcal{E}(\hat{R}_k,\hat{\Omega}_k,P_k)$,
$x_k\in\mathcal{E}_{\Re^6}(0,P_k)$ by the definition of the
uncertainty ellipsoid given in \refeqn{ueso}. This implies that
$A_k^f x_k$ lies in the following uncertainty ellipsoid
\begin{align*}
A_k^f x_k&\in\mathcal{E}_{\Re^6}\!\parenth{0,A_k^f P_k
\parenth{A_k^f}^T}.
\end{align*}
Thus, the uncertainty matrix at the $(k+1)$th step is given by
\begin{align}
P_{k+1}^f & = A_k^f P_k \parenth{A_k^f}^T.\label{eqn:Pkpf}
\end{align}
The above equation is then applied repeatedly to find the
uncertainty matrix at the $(k+l)$th step.

In summary, the uncertainty ellipsoid at the $(k+l)$th step is
computed using \refeqn{findf}, \refeqn{updateR}, \refeqn{updatePi},
and \refeqn{Pkpf} as:
\begin{align}\label{eqn:flow}
    (R_{k+l},\Omega_{k+l})\in\mathcal{E}(\hat{R}_{k+l}^{f},\hat{\Omega}_{k+l}^{f},P_{k+l}^f),
\end{align}

\subsection{Measurement update}
The measurement update finds an uncertainty ellipsoid in the state
space using the measurement and the measurement error models
described in Section \ref{sec:mea}. A feasible set of rotation
matrices that is compatible with the single direction measurement is
described in \refeqn{Ractm}. We find an expression for the
measurement uncertainty ellipsoid such that it contains the set
described by \refeqn{Ractm}.

Elements in the measurement uncertainty ellipsoid are expressed as
\begin{align}\label{eqn:Rmea}
    R_{k+1} = \hat R^m_{k+l} \exp S(\zeta^m_{k+l})
\end{align}
for the center $\hat R^m_{k+l}\in\SO$ and some
$\zeta^m_{k+l}\in\Re^3$. We omit the subscript $(k+l)$ hereafter for
convenience, and it is assumed that the direction is measured at the
$(k+l)$th step.

\textit{Center:} Comparing \refeqn{Ractm} and \refeqn{Rmea}, we
choose the center of the measurement uncertainty ellipsoid as
\begin{align}
    \hat R^m &= \tilde R^\circ,\nonumber\\
    &=\exp\!\bracket{\cos^{-1}(\tilde b^Te)\,\, S\!\parenth{\frac{\tilde b\times
    e}{\|\tilde b\times e\|}}}\exp S(\theta^\circ \tilde b),\label{eqn:Rhatm}
\end{align}
for the constant $\theta^\circ\in S^1$ which is determined by the
following filtering procedure.

\textit{Uncertainty Matrix:} From \refeqn{Ractm} and \refeqn{Rmea},
we have
\begin{align*}
    \exp\! S(\zeta^m)\exp\!\!\bracket{-S(c\tilde b - \nu)}=\exp\!\!\bracket{\theta S\parenth{(I_{3\times
    3}+S(\nu))\tilde b}}.
\end{align*}
Since the vectors $\zeta^m$ and $\zeta^\circ=c\tilde b - \nu$ are
assumed to be small, the above equation is approximated as
\begin{align*}
    \zeta^m -(c\tilde b -\nu)= \theta\tilde b + \theta S(\nu)\tilde
    b,
\end{align*}
which can be rewritten as
\begin{align*}
    \zeta^m-(\theta+c)\tilde b & = -\nu -\theta S(\tilde b)\nu,\\
    & = (\theta-1)\nu-\theta (I_{3\times 3}+S(\tilde b))\nu.
\end{align*}
Since $\nu\in\mathcal{E}_{\Re^3}(0_3,S)$ and $\theta\in\S^1$, the
terms in the right hand side satisfy
\begin{gather*}
    (\theta-1)\nu \in \mathcal{E}_{\Re^3}(0_3,(1+\pi)^2 S),\\
\theta (I_{3\times 3}+S(\tilde b))\nu \in
\mathcal{E}_{\Re^3}(0_3,\pi^2 \mathcal{A}^{m,T}S\mathcal{A}^{m,T}),
\end{gather*}
where $\mathcal{A}^{m}=I_{3\times 3}+S(\tilde b)\in\Re^{3\times 3}$.
Therefore, the vector $\zeta^m-(\theta+c)\tilde b$ lies in an
ellipsoid containing the vector sum of the above two ellipsoids. The
expressions for the minimal ellipsoid containing the vector sum of
two ellipsoids are given in \cite{jo:MaNo1996}. Using the results,
we have
\begin{align}\label{eqn:zetam0}
    \zeta^m-(\theta+c)\tilde b \in \mathcal{E}_{\Re^3}(0_3,P^m_0),
\end{align}
where
\begin{gather*}
    P^m_0=(1+q^{-1})Q^{1}+(1+q)Q^{2},\quad
    q=\sqrt{\frac{\tr{Q^{1}}}{\tr{Q^{2}}}},\\
    Q^{1}=(1+\pi)^2 S,\quad
    Q^{2}=\pi^2 \mathcal{A}^{m,T}S\mathcal{A}^{m}.
\end{gather*}

From \refeqn{zetam0}, we can guarantee that the vector $\zeta^m$
lies in an ellipsoid containing the following union of the sets
\begin{align*}
    \mathcal{E}_{\Re^3}(-\pi\tilde b,P^m_0)\cup\mathcal{E}_{\Re^3}(\pi\tilde
    b,P^m_0).
\end{align*}
This is a consequence of the fact that an ellipsoid is convex and
the assumption $\abs{c}\ll \pi$. The ellipsoid that contains the
union of two ellipsoids is obtained numerically by the LMI approach
presented in~\cite{BoyGhaFer.BK94}.
\begin{align}\label{eqn:Pm}
    \mathcal{E}_{\Re^3}(0_3,P^m)\supset\parenth{\mathcal{E}_{\Re^3}(-\pi\tilde b,P^m_0)\bigcup\mathcal{E}_{\Re^3}(\pi\tilde
    b,P^m_0)}.
\end{align}

In summary, a single direction measurement with small error
guarantees that the rotation matrix is expressed as \refeqn{Rmea},
where the center $\hat R^m$ is given by \refeqn{Rhatm}, and the
vector $\zeta^m$ lies in the uncertainty ellipsoid given by
\refeqn{Pm}.

\subsection{Filtering procedure}
The filtering procedure finds a new uncertainty ellipsoid compatible
with both the predicted uncertainty ellipsoid and the measured
uncertainty ellipsoid. The intersection of two ellipsoids is
generally not an ellipsoid. We find a minimal uncertainty ellipsoid
containing the intersection.

The predicted uncertainty ellipsoid is based on $\hat R^f$ and the
measurement ellipsoid is based on $\hat R^m$. In the following
development, we assume that the difference between the rotation
matrices $\hat R^f$ and $\hat R^m$ is small. Here we find a value of
$\theta^\circ\in\S^1$ at \refeqn{Rhatm} such that the difference is
minimized. Define an index $\mathcal{J}=\tr{I_{3\times 3}-\hat
R^{f,T} \hat R^m}$. A standard variational approach with the use of
Rodriguez formula shows that the index is minimized when
\begin{gather*}
    \theta^\circ=-\tan^{-1}\frac{\tr{\hat R^{f,T}\tilde R^\vartriangle S(\tilde b)}}{\tr{\hat R^{f,T}\tilde R^\vartriangle S(\tilde b)^2}},\\
     \tr{\hat R^{f,T}\tilde R^\vartriangle S(\tilde
    b)}\sin\theta^\circ-\tr{\hat R^{f,T}\tilde R^\vartriangle S(\tilde
    b)^2}\cos\theta^\circ >0,
\end{gather*}
where $\tilde R^\vartriangle\in\SO$ is the first exponential of
\refeqn{Rhatm}. The first equation is obtained by the optimality
condition $\deriv{\mathcal{J}}{\theta^\circ}=0$, and the second
inequality is obtained by
$\frac{\partial^2\mathcal{J}}{\partial(\theta^\circ)^2}>0$. These
conditions define the value of $\theta^0\in\S^1$ uniquely.

We find a minimal ellipsoid containing the intersection of the
predicted uncertainty ellipsoid and the measurement uncertainty
ellipsoid. An element in the predicted uncertainty ellipsoid,
$(R^f,\Omega^f)\in\mathcal{E}(\hat{R}^{f},\hat{\Omega}^{f},P^f)$,
can be written as
\begin{align}
R^f & = \hat{R}^{f} e^{S(\zeta^f)},\label{eqn:Rfm}\\
\Omega^f & = \hat{\Omega}^{f} + \delta{\Omega}^f,\label{eqn:Pifm}
\end{align}
for some $(\zeta^f,\delta\Omega^f)\in\mathcal{E}_{\Re^6}(0_{6\times
1},P^f)$. We find an equivalent expression based on the measurement
ellipsoid center $\hat R^m$. Define $\hat\zeta^{mf}\in\Re^3$ such
that
\begin{align}
\hat{R}^f&=\hat{R}^m e^{S(\hat\zeta^{mf})}.\label{eqn:zetamfm}
\end{align}
Thus, $\hat\zeta^{mf}$ represents the difference between the centers
of the two ellipsoids. Substituting \refeqn{zetamfm} into
\refeqn{Rfm},
\begin{align*}
R^f & = \hat{R}^m e^{S(\hat\zeta^{mf})} e^{S(\zeta^f)},\nonumber\\
& \simeq \hat{R}^m e^{S(\hat\zeta^{mf}+\zeta^f)},
\end{align*}
where we assumed that $\hat\zeta^{mf}, \zeta^f$ are sufficiently
small to obtain the second equality. Thus, the uncertainty ellipsoid
obtained by the flow update,
$\mathcal{E}(\hat{R}^{f},\hat{\Omega}^{f},P^f)$ is identified by the
center $(\hat{R}^m,\hat{\Omega}^f)$ and the uncertainty ellipsoid in
$\Re^6$.
\begin{align*}
(\zeta^{mf},\delta\Omega^{f})\in\mathcal{E}_{\Re^6}(\hat{x}^{mf},P^f),
\end{align*}
where $\hat{x}^{mf}=[\hat\zeta^{mf};0_{3}]$, and
$S(\zeta^{mf})=\mathrm{logm} \parenth{\hat{R}^{m,T} R^f}\in\so$,
$\delta\Omega^{f}=\Omega^f-\hat{\Omega}^{f}\in\Re^3$.

We seek a minimal ellipsoid that contains the intersection of the
following uncertainty ellipsoids.
\begin{align}
\parenth{\mathcal{E}_{\Re^3}(0_3,P^m)\bigcap
\mathcal{E}_{\Re^6}(\hat{x}^{mf}
,P^f)}\subset\mathcal{E}_{\Re^6}(\hat{x},P),
\end{align}
where $\hat{x}=[\hat\zeta^T,\delta\hat\Omega^T]^T\in\Re^6$.
Expressions for a minimal ellipsoid containing the intersection of
two ellipsoids are presented in~\cite{jo:MaNo1996}. Using those
results, $\hat{x}$ and $P$ are given by
\begin{align*}
\hat{x}&=(I_{6\times 6}-LH)\hat{x}^{mf},\\
P&=\beta(r) [(I-LH)P^f(I-LH)^T+r^{-1}LP^mL^T],
\end{align*}
where $L\in\Re^{6\times 3}$, $H\in\Re^{3\times 6}$, and
$\beta(r)\in\Re$ are given by
\begin{align*}
L & = P^{f} H^T[HP^{f}H^T + r^{-1} P^m]^{-1},\\
H & = [I_{3\times 3}, 0_{3\times 3}]^T,\\
\beta(r) & = 1 + r - (\hat{x}^{mf})^T H^T [HP^fH^T + r^{-1}P^m]^{-1}
H\hat{x}^{mf},
\end{align*}
for a constant $r$, which is chosen such that $\tr{P}$ is minimized.

In summary, a new uncertainty ellipsoid at the $(k+l)$th step is
given by
\begin{align}
    (R_{k+1},\Omega_{k+1})\in\mathcal{E}(\hat{R}_{k+1},\hat{\Omega}_{k+1},P_{k+1}),
\end{align}
where
\begin{align}
\hat{R}_{k+1}&=\hat{R}_{k+1}^m e^{S(\hat\zeta)},\\
\hat{\Omega}_{k+1}& = \hat{\Omega}_{k+1}^f + \delta\hat\Omega,\\
P_{k+1}& = P.
\end{align}
The entire procedure  is repeated whenever a new measurement is
available.

The steps outlined above define a dynamic filter. 
The center of the uncertainty ellipsoid is considered as a point
estimate of the attitude and the angular velocity at the $(k+l)$th
step. The uncertainty matrix represents the characteristics of the
uncertainty, and the size of the uncertainty matrix represents the
accuracy of the estimate. If the size of the uncertainty ellipsoid
is small, we conclude that the estimate is accurate. This estimation
is optimal in the sense that the size of the filtered uncertainty
ellipsoid is minimized.

\subsection{Properties of the estimator}
The notable feature of this attitude estimator is that it requires a
single direction measurement. Current attitude estimators based on
the solution of Wahba's problem require at least two direction
measurements at each instant. A single direction measurement
provides only a two-dimensional constraint for the six-dimensional
tangent bundle. The information obtained from the attitude dynamics
is utilized, together with the measurement, in order to estimate the
attitude and the angular velocity of the rigid body. In this paper,
it is assumed that the angular velocity is not measured, but the
current results can be readily extended to incorporate angular
velocity measurements.

This attitude estimator has no singularities since the attitude is
represented by a rotation matrix, and the geometric structure of the
rotation matrix is preserved since it is updated by the
structure-preserving Lie group variational integrator. The presented
estimator can be used for highly nonlinear large angle maneuvers of
a rigid body. It is also robust to the distribution of the
measurement noise since we only use ellipsoidal bounds on the noise.
The measurements need not be periodic, the estimation is repeated
whenever new measurements become available.

\section{Numerical Example}
Numerical simulation results are given for the estimation of the
attitude dynamics of an uncontrolled rigid spacecraft in a circular
orbit about a large central body, including gravity gradient
effects. The details of the on orbit spacecraft model are presented
in~\cite{pro:acc06}.

The mass, length and time dimensions are normalized by the mass of
the spacecraft, the maximum length of the spacecraft, and the
orbital angular velocity, respectively. The moment of inertia of the
spacecraft is chosen as ${J}=\mathrm{diag}\bracket{1,\,2.8,\,2}$.
The maneuver is a large attitude change completed in a quarter of
the orbit. The initial conditions are chosen as
\begin{alignat*}{2}
R_0 & = \begin{bmatrix}0.707 & -0.707 & 0\\
                        0.707 & 0.707 & 0\\
                        0 & 0 & 1\end{bmatrix},&
\Omega_0&=[2.32,\,0.45,\,-0.59],\\
\hat R_0 & = I_{3\times 3},&\quad
\hat{\Omega}_0&=[2.12,\,0.55,\,-0.89].
\end{alignat*}
The corresponding initial estimation errors are
$\norm{\zeta_0}=45\,\mathrm{deg}$,
$\norm{\delta\Omega_0}=21.43\frac{\pi}{180}\,\mathrm{rad/s}$. Note
that the actual initial attitude is opposite to the estimated
initial attitude. The initial uncertainty matrix is given by
\begin{align*}
P_0 = \mathrm{diag}\bracket{2.28,2.28,2.28,0.82,0.82,0.82},
\end{align*}
so that $x_0^TP_0^{-1}x_0=0.857\leq 1$.

We assume that the measurements are available twenty times. The
inertial direction $e$ to a known point is chosen from the columns
of the following matrix.
\begin{align*}
E=\begin{bmatrix}1 & 0 & 0 & 0.7071 & -0.7071 & 0.5 & -0.5\\
    0 & 1 & 0& 0.7071 & 0.7071 & 0.5 & 0.5\\
    0 & 0& 1 & 0 & 0 & 0.7071 & 0.7071\end{bmatrix}.
\end{align*}
A simple adaptive scheme is developed to choose the best inertial
direction as the spacecraft rotates. The uncertainty matrix for the
measurement noise is given by
\begin{align*}
S_k=\parenth{0.2\frac{\pi}{180}}^2.
\end{align*}
The direction measurement noise is normally distributed in the
simulation.

\begin{figure}
    \centerline{\subfigure[Estimation error $\norm{\zeta_k}$, $\norm{\delta\Omega_k}$]{
    \includegraphics[width=0.46\columnwidth]{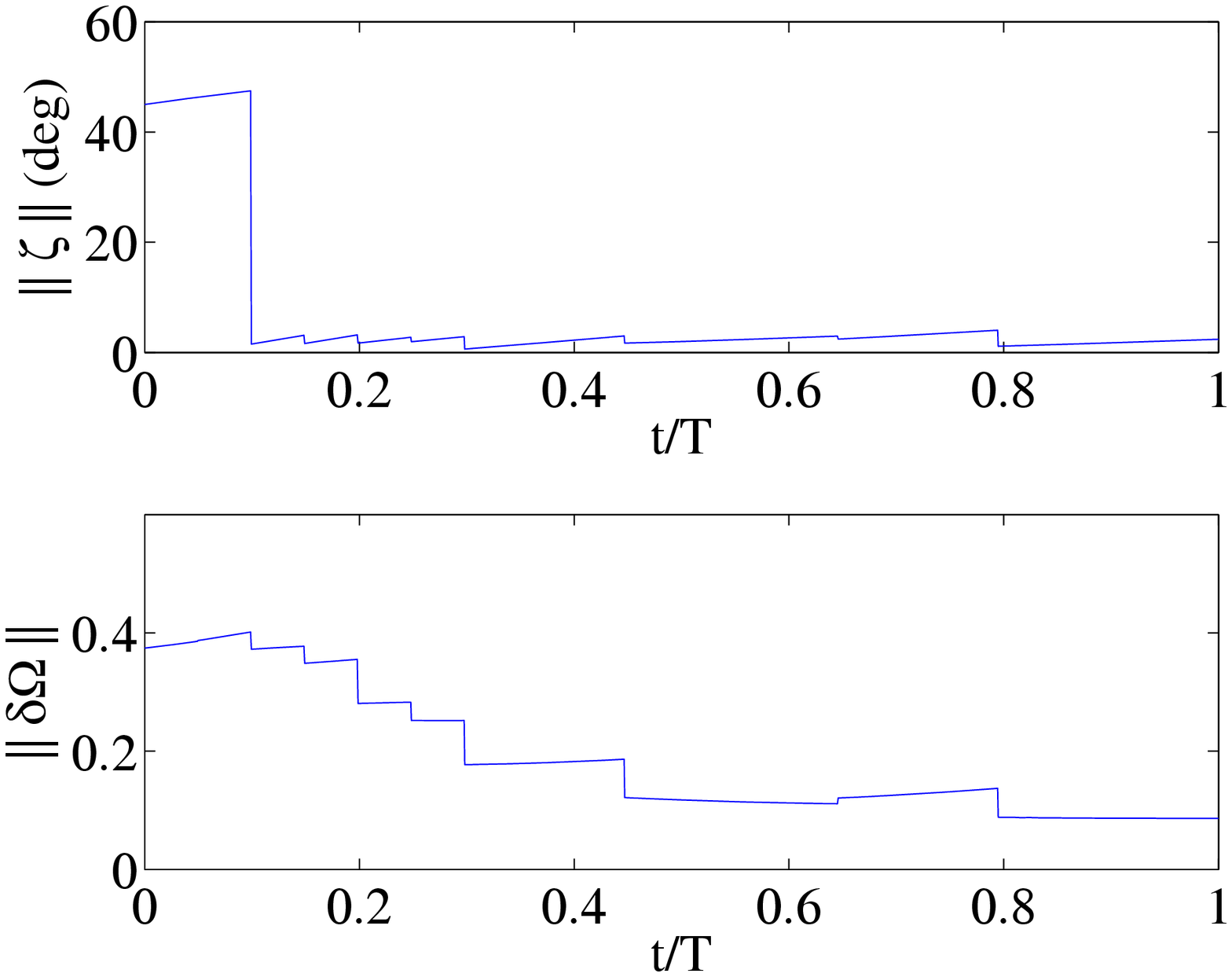}}
    \hfill
    \subfigure[Size of uncertainty $\tr{P_k}$]{
    \includegraphics[width=0.45\columnwidth]{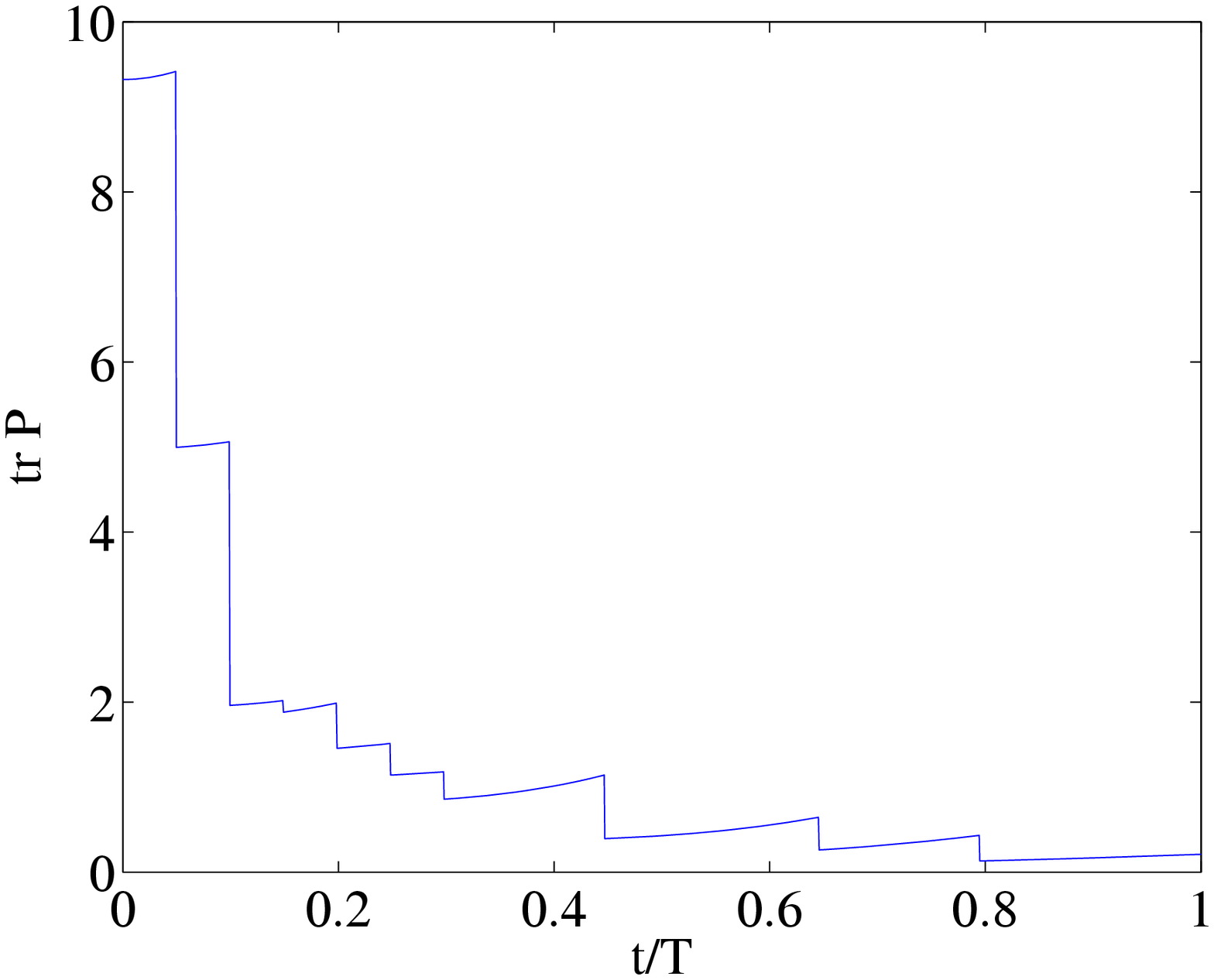}}
    }
    \caption{Estimation errors and uncertainty}\label{fig:sim}
\end{figure}

\reffig{sim} shows simulation results, where the left figure shows
the attitude estimation error and the angular velocity estimation
error, and the right figure shows the size of the uncertainty
ellipsoid. The estimation errors and the size of the uncertainty are
reduced rapidly after the first few measurements; the estimation
error for the angular velocity converges relatively slowly since the
angular velocity is not measured directly. The terminal attitude
error, and the terminal angular velocity error are less than
$2.3\,\mathrm{deg}$, and $0.08\,\mathrm{rad/s}$, respectively.

\section{Conclusions}
A deterministic attitude estimator for a rigid body under an
attitude dependent potential is developed. This estimator requires
only a single direction measurement to a known reference point at
each measurement instant. A feasible set of rotation matrices
compatible with the measurement is described in terms of Lie algebra
elements, and it is compared with an uncertainty ellipsoid obtained
from an attitude dynamics model, in order to obtain an updated
attitude estimate. The attitude is globally represented by a
rotation matrix, and the geometric structure of the rotation matrix
is preserved by using a Lie group variational integrator.

\bibliography{ue}
\bibliographystyle{IEEEtran}

\end{document}